\documentclass{article}

\usepackage{color}
\usepackage{amssymb}
\usepackage{amsmath}

\newtheorem{theor}{Theorem}
\newtheorem{prop}{Proposition}
\newtheorem{lm}{Lemma}

\begin{document}

\title{On some Godbillon-Vey classes\\
of a family of regular foliations}
\author{ Cristian Ida and Paul Popescu}
\date{}
\maketitle

\begin{abstract}
The aim of the paper is to construct some Godbillon-Vey classes of a family
of regular foliations, defined in the paper. These classes are cohomology
classes on the manifold or on suitable open subsets. Some examples are also
considered.
\end{abstract}



\section{Introduction}

The families of regular foliations considered in the paper are regular
foliations on open subsets such that all the induced leaves on an
intersection set give a system of subfoliations as in \cite{Do, Ho} (i.e.
the induced larger-sized leaves are saturated with smaller-sized ones; see
conditions (F1)--(F3) in the next section). The resulting geometric
distribution, given by the tangent subspaces to leaves of maximal dimension,
is a singular one (1. of Proposition \ref{prSuSt}). Assuming that any
intersection is saturated by whole leaves, particular classes of
Stefan-Sussmann foliations are obtained (2. of Proposition \ref{prSuSt}),
called here singular foliations that are locally regular.

A tool used to extend Godbillon-Vey forms, on a stratum with a non-minimal
dimensional leaves, is the existence of a basic test function on the
complement of the stratum. We call a test function, according to a closed
subset $M_{0}\subset M$, a smooth real function that has $M_{0}$ as its set
of zeros. The existence of a general test function follows from a classical
results of Whitney and some properties of extension of smooth sections on
closed subsets (see \cite{Lee, Mu, Wi}, but in a slight different form).
Using the line of \cite[Section 4]{DLPR}, we give a proof in Proposition \ref%
{prtest}.

The main constructions in the paper are performed in the fourth section. The
most important one is that of the Godbillon-Vey class of leaves of minimal
dimension in $M$ and in $\Sigma _{\geq r_{i}}$ (Theorem \ref{thmain}), where
we prove that the Godbillon-Vey form of the leaves extends to a global
cohomology class $GV_{\min }(\mathcal{F})\in H^{1+2q_{\max }}(M)$ (for the
leaves of minimal dimension on $U_{0}$) and to some Godbillon-Vey classes $%
GV_{\min }(\mathcal{F}_{\Sigma _{\geq r_{i}}})\in H^{2\left( m-r_{i}\right)
+1}(\Sigma _{\geq r_{i}})$ (for the other leaves on $U_{i}$, $i>0$). In the
case when there is a basic test function of $M\backslash U_{i}$, then one
get a cohomology class on $M$ (Proposition \ref{prbaz}).

Two cases are considered in the last section. First, given a regular
foliation $\mathcal{F}_{0}$ on $M$, one can easily construct a family of
regular foliations on $M$ (for example, adding in a suiatble open set a
trivial foliation with one leaf), such that its Godbillon-Vey class $%
GV_{\min }(\mathcal{F})$ is the same as $GV(\mathcal{F}_{0})$, the usual
Godbillon-Vey class of $\mathcal{F}_{0}$ (Proposition \ref{prex01}). Thus if
the the Godbillon-Vey class of $\mathcal{F}_{0}$ is non-trivial, also is
that of the family $\mathcal{F}$. Second, we prove that if $0$ is a regular
value for the (weak) test basic function $\varphi _{i}$, then the cohomology
class $[\overline{\nu }_{i}]\in H^{2q_{i}+1}(M)$ vanishes (Proposition \ref%
{prex02}).

Looking at the first example, it seems likely to find a non-trivial family\
of regular foliations, maybe a singular foliation, that is locally regular,
having a more complicated structure and a non-trivial Godbillon-Vey class.
The second example shows that a non-trivial Godbillon-Vey class can be found
not for a regular (weak) test function, possible for a strong one. We let it
as an open problem.  

\section{Families\ of regular foliations}

Let $M$ be a differentiable manifold. Let us suppose that there is an open
cover $\{U_{i}\}_{i\in I}$ of $M$ such that the following three conditions
hold:

(F1) -- on every $U_{i}$ there is a regular foliation $\mathcal{F}_{i}$
having $r_{i}$ as dimension of leaves,

(F2) -- if $i\neq j$ then $r_{i}\neq r_{j}$ and

(F3) -- if $U_{i}\cap U_{j}\neq \emptyset $, $r_{i}<r_{j}$, then $U_{i}\cap
U_{j}$ is saturated by open subsets of leaves of $\mathcal{F}_{j}$ and every
such open set is saturated to its turn by open subsets of leaves of $%
\mathcal{F}_{i}$.

We can consider a stronger condition than (F3) as:

(F3') -- if $U_{i}\cap U_{j}\neq \emptyset $, $r_{i}<r_{j}$, then $U_{i}\cap
U_{j}$ is saturated by leaves of $\mathcal{F}_{j}$ and every such leaf of $%
\mathcal{F}_{j}$ is saturated to its turn by leaves of $\mathcal{F}_{i}$.

It is easy to see that $I$ is a finite set, $I=\overline{0,k}$. The rank of
a point $x\in M$ is $r(x)=\max \{r_{i}:x\in U_{i}\}$; if $r(x)=r_{i}$, then
and we denote by $\mathcal{D}_{x}$ the tangent space to the leaf of $%
\mathcal{F}_{i}$. We denote by $\mathcal{R}=\{r(x)=\dim \mathcal{D}_{x}:x\in
M\}$. If $S\subset M$ then $\mathcal{D}_{S}=\underset{x\in S}{\cup }{}%
\mathcal{D}_{x}$ denotes the restriction of $\mathcal{D}$ to $S$. Let $\mathcal{R}=$ 
$\{r_{i}\}_{i=\overline{0,k}}$, where $r_{\min }=r_{0}<r_{1}<\cdots
<r_{k}=r_{\max }$. For $r_{i}\in \mathcal{R}$, we denote by $\Sigma
_{r_{i}}=\{x\in M:$ $\dim \mathcal{D}_{x}=r_{i}\}$, $\Sigma _{<r_{i}}=\{x\in
M:$ $\dim \mathcal{D}_{x}<r_{i}\}$, $\Sigma _{\leq r_{i}}=\{x\in M:$ $\dim 
\mathcal{D}_{x}\leq r_{i}\}=\Sigma _{r_{i}}\cup \Sigma _{<r_{i}}$, $\Sigma
_{>r_{i}}=\{x\in M:$ $\dim \mathcal{D}_{x}>r_{i}\}$, $\Sigma _{\geq
r_{i}}=\{x\in M:$ $\dim \mathcal{D}_{x}\geq r_{i}\}=\Sigma _{r_{i}}\cup
\Sigma _{>r_{i}}$. We say that the subset $\Sigma _{r_{\min }}$ is the \emph{%
minimal set} and $\Sigma _{r_{\max }}$ is the maximal set. The subsets $%
\Sigma _{<r_{i}}$ and $\Sigma _{\leq r_{i}}$ are closed subsets and their
complements, the sets $\Sigma _{\geq r_{i}}$ and $\Sigma _{>r_{i}}$ are open
closed subsets in $M$. The subset $\Sigma _{r_{i}}\subset \Sigma _{\geq
r_{i}}$ is the minimal subset of $\mathcal{D}_{|\Sigma _{\geq r_{i}}}$ and $%
\Sigma _{>r_{i}}$ is void if $i=k$ and is equal to $\Sigma _{\geq r_{i+1}}$
if $0\leq i<k$. We say also that the leaves of $\mathcal{F}_{i}$ are \emph{%
leaves of minimal dimension}.

The assignment of a vector subspace $\mathcal{D}_{x}\subset T_{x}M$, $%
(\forall )x\in M$, gives a \emph{singular distribution }$\mathcal{D}$ on $M$%
, $\mathcal{D}=\underset{x\in M}{\cup }{}\mathcal{D}_{x}\subset TM$. We
denote by $\Gamma _{loc}(\mathcal{D})$ the set of local smooth vector fields
tangent to $\mathcal{D}$ in every point where they are defined. One say that 
$\mathcal{D}$ is:

-- \emph{smooth,} if $\mathcal{D}_{x}$ is spanned by some restrictions to $x$
of some smooth local vector fields from $\Gamma _{loc}(\mathcal{D})$, $%
(\forall )x\in M$;

-- (\emph{completely}) \emph{integrable}, if $\mathcal{D}$ is smooth and
there is a partition of $M$ in immersed submanifolds $L\subset M$ such that
if $x\in L$, then $\mathcal{D}_{x}=T_{x}L$.

(See, for example \cite{Da, Va1} for more details.)

\begin{prop}
\label{prSuSt}.

\begin{enumerate}
\item Assuming the conditions (F1), (F2) and (F3), then $\mathcal{D}$ is a
smooth singular distribution on $M$.

\item Assuming the conditions (F1), (F2) and (F3'), then the singular
distribution $\mathcal{D}$ is integrable.
\end{enumerate}
\end{prop}

\emph{Proof.} Let $x\in M$ and a regular foliate chart of the leaf $F_{i}$
of $\mathcal{F}_{i}$ that contain $x$, where $r(x)=r_{i}$. The condition
(F3) implies that the canonical tangent vectors to $F_{i}$ belong to $\Gamma
_{loc}(\mathcal{D})$ and their restrictions to $x$ generate $T_{x}F_{i}=%
\mathcal{D}_{x}$. Assuming supplementary the condition (F3'), then this
local chart is also one corresponding to a singular Stefan-Sussmann
foliation on $M$ (according for example to \cite{Va1}) that is tangent to $%
\mathcal{D}$. $\Box $

We say that

-- the conditions (F1), (F2) and (F3) define a \emph{family\ of regular
foliations} and

-- the conditions (F1), (F2) and (F3') define a \emph{singular foliation
that is locally regular}.

For a family of regular foliations, we can define the \emph{leaf} of $x\in M$
as the leaf $F_{i}$ of $\mathcal{F}_{i}$ that contains $x$, of maximal
dimension $r(x)=r_{i}$. Moreover, in general a non-ambiguous leaf can be
defined only for totally integrable foliations.

Notice that the conditions (F1), (F2) and (F3) does not always assure that $%
\mathcal{D}$ (defined as above) is integrable. Indeed, consider the open
cover of $I\!\!R^{2}$ given by $U_{1}=\{(x,y)\in I\!\!R^{2},$ $x>0\}$ and $%
U_{2}=\{(x,y)\in I\!\!R^{2},$ $x<1\}$. Let us consider the foliation $%
\mathcal{F}_{1}$ by one leaf on $U_{1}$ and the foliation $\mathcal{F}_{2}$
by horizontal lines $y=const.$ on $U_{2}$. The conditions (F1)-(F3) ale
fulfilled, but the condition (F3') is not fulfilled. It generates a singular
smooth distribution $\mathcal{D}$ that is not integrable, generated by the
vector fields $X_{1}=\dfrac{\partial }{\partial x}$ and $X_{2}=\varphi (x)%
\dfrac{\partial }{\partial y}$, where $\varphi $ vanishes for $x\leq 0$ and $%
\varphi (x)=e^{-\frac{1}{x}}$ for $x>0$.

Let us consider some other examples.

-- Given a family\ of regular foliations (or a singular foliation that is
locally regular), the open set $\Sigma _{\geq r}$ is saturated by leaves of $%
\mathcal{F}_{i}$, where $r_{i}\geq r$, thus a family\ of regular foliations
(or a singular foliation that is locally regular) $\mathcal{F}_{\geq r}$ is
induced. In particular $\mathcal{F}_{\geq r_{k}}=\mathcal{F}_{r_{k}}$ on $%
\Sigma _{\geq r_{k}}=\Sigma _{r_{\max }}$ is regular.

-- A regular foliation on $M$ is a singular foliation that is locally
regular, when all the points have the same rank, equal to the dimension of
the leaves (i.e. of the foliation).

-- A non-trivial example is given by the foliation of $I\!\!R^{n}$ by concentric spheres
(as leaves of dimension $n-1$) and the origin (as a leaf of dimension $0$) that
is a singular foliation that is locally regular. An other non-trivial
example is a singular foliation having as leaves concentric spheres, as in
the previous example (of dimension $n-1$), outside a compact ball $\bar{B}(%
\bar{0},\rho )\in I\!\!R^{n}$, $\rho >0$, while $\bar{B}(\bar{0},\rho )$ is
a union of points (as leaves of dimension $0$).

-- A singular Stefan-Sussmann foliation on $M$ that has $\mathcal{R}=\{0,r\}$%
, where $0<r\leq m=\dim M$ is locally regular. In general, consider a
regular foliation $\mathcal{F}_{U}$ on an open subset $U\subset M$, such
that the dimension of leaves is $r$, where $0<r\leq m$. The partition of $M$
by the leaves of $U$ and by the points of $\Sigma _{0}=M\backslash U$ gives
a locally regular Stefan-Sussmann foliation on $M$. The singular
distribution has $\mathcal{R}=\{0,r\}$. Notice that any singular
Stefan-Sussmann foliation on $M$ that has $\mathcal{R}=\{0,r\}$ can be
obtained in this way.

-- Consider a regular foliation $\mathcal{F}_{U}$ on an open subset $%
U\subset M$, such that the dimension of leaves is $r$, where $0\leq r<m$.
Let $\Sigma _{0}\subset U$ be a closed subset of $M$, saturated or not by
leaves of $\mathcal{F}_{U}$. The partition of $M$ by the leaves of $\mathcal{%
F}_{\Sigma _{0}}$ and the leaf $\Sigma _{1}=M\backslash \Sigma _{0}$ gives a
family\ of regular foliations This is a singular foliation that is locally
regular only if $\Sigma _{0}$ is saturated by the leaves of $\mathcal{F}_{U}$%
, when it gives a locally regular Stefan-Sussmann foliation on $M$. This
singular distribution has $\mathcal{R}=\{r,m\}$.

-- Consider some open subsets $U_{1}$, $U_{2}\subset M$ and a regular
foliation $\mathcal{F}_{1}$ on $U_{1}$; we suppose that $U_{1}\cap U_{2}\neq
\emptyset $ and $U_{1}\cup U_{2}\neq M$. Denote by $\Sigma _{0}=M\backslash
(U_{1}\cup U_{2})$ and let $U_{0}\supset \Sigma _{0}$ be an open set. We
consider on $U_{0}$ and $U_{2}$ the trivial foliations $\mathcal{F}_{0}$ and 
$\mathcal{F}_{2}$ respectively, where $\mathcal{F}_{0}$ has points as leaves
and $\mathcal{F}_{2}$ has one leaf. It follows a family\ of regular
foliations. If $U_{1}\cap U_{2}$ is saturated by leaves of $\mathcal{F}_{1}$%
, then the family\ of regular foliations is a singular foliation that is
locally regular.

The suspension constructed for regular foliations (as, for example, in \cite[%
2.7, 2.8]{Go}) can be extended to a family\ of regular foliations, as
follows. Let $B$ and $M$ be two manifolds and $\mathcal{F}$ be a family of
regular foliations or a singular foliation that is locally regular. Let us
suppose that $\rho :\pi _{1}(B)\rightarrow Diff(M)$ is a representation
(i.e. a group morphism) such that every diffeomorphism $\rho (g)\in Diff(M)$
invariate an open neighborhood $U_{k}$ of $\Sigma _{k}$, as well as the
leaves of the foliation $\mathcal{F}_{k}$ on $U_{k}$ that restricts to the
leaves on $\Sigma _{k}$. If we denote by $\tilde{B}$ the universal simple
connected cover of $B$, then the suspension space is the quotient space $S=(
\tilde{B}\times M)/^{\sim}$ of the equivalence relation $(\tilde{b},m)%
\sim(\tilde{b}g,\rho (g)^{-1}m)$, $g\in \pi _{1}(B)$, on $\tilde{B}%
\times M$. As in the classical case, one can first consider on $\tilde{B}%
\times M$ the product foliations $\mathcal{F}_{0}$ of the foliation by one
leaf on $\tilde{B}$ and the foliations $\mathcal{F}_{i}$ on $M$. A family
of regular foliations or a singular foliation that is locally regular
(accordingly to that on $M$) is induced on the quotient space $S$; the
leaves, the sets $\Sigma _{k^{\prime }}$ of the leaves of a same dimension $%
k^{\prime }$ and the open neighborhoods $U_{k^{\prime }}^{\prime }$ of $%
\Sigma _{k^{\prime }}$ are naturally induced.

As a particular case, consider an open subset $U\subset M$, a regular
foliation $\mathcal{F}_{U}$ on $U$ and $f\in Diff(M)$ such that $f(U)=U$ and 
$f$ invariates $\mathcal{F}_{U}$. We can consider an open neighborhood $W$
of the closed set $M\backslash U$ (for example $W=M$) and the trivial
foliation $\mathcal{F}_{W}$ by points on $W$. The leaves of $\mathcal{F}_{U}$
and the points of $M\backslash U$ as $0$-dimensional leaves give a locally
regular Stefan-Sussmann foliation on $M$. The suspension of $f$ is
considered for $B=S^{1}$, $\tilde{B}=I\!\!R$, $\pi _{1}(S^{1})=Z\!\!\!Z$ and
the actions $I\!\!R\times Z\!\!\!Z\rightarrow I\!\!R$, $(x,n)\rightarrow x-n$
and $Z\!\!\!Z\times M\rightarrow M$, $\left( n,m\right) \rightarrow f^{n}(m)$%
.

For example, consider the natural central symmetry $\sigma :S^{n}\rightarrow
S^{n}\subset I\!\!R^{n+1}$, $\sigma (\bar{x})=-\bar{x}$. Consider also two
open spherical caps $C_{1}\subset C_{2}$ centred in the same point $A$ of
the sphere $S^{n}$ and let $C_{1}^{\prime }=\sigma (C_{1})\subset
C_{2}^{\prime }=\sigma (C_{2})$ the symmetric spherical caps centred in $%
A^{\prime }=\sigma (A)$, such that $C_{2}\cap C_{2}^{\prime }\neq \emptyset $%
. Denote by $U_{1}=S^{n}\backslash (\bar{C}_{1}\cup \bar{C}_{1}^{\prime })$
and by $U_{2}=C_{2}\cup C_{2}^{\prime }$. Consider the trivial foliation $%
\mathcal{F}_{2}$ on $U_{2}$ by points and a $k$-regular foliation $\mathcal{F%
}_{1}$ on $U_{1}$ obtained by intersection of $U_{1}$ by $k+1$--parallel
planes that can be parallel or not with the support $n$--hyperplanes of the
spherical caps. Obviously the open sets $U_{1}$ and $U_{2}$, as well as the
foliations $\mathcal{F}_{1}$ and $\mathcal{F}_{2}$ are invariated by $\sigma 
$. One can consider a quotient locally regular foliation on $I\!\!RP^{n}$,
as well as a suspension locally regular foliation on $S=(I\!\!R\times S^{n})/^\sim$ , given by the $Z\!\!\!Z$--action $n\cdot (\alpha ,\bar{x}%
)=(\alpha -n,\sigma ^{n}(\bar{x}))$.

\section{Test functions}

We consider now test functions, that allow us to extend smooth functions and
vector fields.

Let $M_{0}\subset M$ be a closed subset. We say that a real function $%
\varphi \in \mathcal{F}(M)$ is a \emph{weak test function} for $M_{0}$ if $%
M_{0}=\varphi ^{-1}(0)$ (i.e. $\varphi (x)=0$ iff $x\in M_{0}$). We say that
a weak test function is a \emph{strong test function} for $M_{0}$ if,
additionally, its values are in $[0,1]$ and all its differentials vanish in
every $x\in M_{0}$. The existence of test functions is an important tool
used in the sequel.

The following simple Lemma shows that the existence of a weak test function
gives a strong one.

\begin{lm}
\label{lmsmooth}Let $\psi _{0}:I\!\!R\rightarrow \lbrack 0,1]$ be smooth function
such that $\psi _{0}(t)=0$ iff $t=0$ and all the derivatives of $\psi _{0}$
vanish in $t=0$. Then for every function $f:M\rightarrow I\!\!R$ the
function $F=\psi _{0}\circ f$ has the same zeros as $f$ and all the
differentials of $F$ vanish in its zeros.
\end{lm}

Notice that a function $\psi _{0}$ as in Lemma \ref{lmsmooth} is 
\begin{equation}
\psi _{0}(t)=\left\{ 
\begin{array}{l}
\frac{e^{-\frac{1}{t^{2}}}}{1+e^{-\frac{1}{t^{2}}}}~\mathrm{if}\ t\neq 0, \\ 
0~\,\,\,\,\,\,\,\,\,\,\,\,\,\,\,\,\mathrm{if}\ t=0.%
\end{array}%
\right. .  \label{ps0}
\end{equation}

A first fact is the existence of a  weak test function $\varphi _{M_{0}}$ for
any closed subset $M_{0}\subset M$, i.e. a positive smooth real function on $%
M$, having the set of zeros exactly $M_{0}$. The existence follows from a
classical results of Whitney and some properties of extension of smooth
sections on closed subsets (see \cite{Lee, Mu, Wi}), but in a slight
different form. We give a proof below, in line of \cite[Section 4]{DLPR}.

\begin{prop}
\label{prtest}Let $M$ be a differentiable manifold and $M_{0}\subset M$ be a
closed subset. Then there is a (weak, strong) test function for $M_{0}$.
\end{prop}

\emph{Proof.} We can proceed as in \cite[Section 4]{DLPR} reducing the
problem to the case when $M=I\!\!R^{n}$ and considering $M$ properly
embedded in $I\!\!R^{k}$ for some $k$. Then $M_{0}\subset I\!\!R^{k}$ is
also closed. A test function on $I\!\!R^{k}$ for $M_{0}$ reduces to $M$ also
to a test function for $M_{0}$. Since $M_{0}$ is a closed set, then $%
M_{1}=I\!\!R^{k}\backslash M_{0}$ is an open subset of $I\!\!R^{k}$. For any
point $p\in M_{1}$ there is a ball $B_{p}=B(p,2r)\subset M_{1}$. We
denote by $B_{p}^{\prime }=B(p,r)$ and we consider a bump function $\psi
_{p}:M\rightarrow \lbrack 0,1]$ such that its support is $\bar{B}_{p}=\bar{B}%
(p,2r)$, its values are $0$ outside $B_{p}$ (i.e. on $I\!\!R^{n}\backslash
B_{r}$), $1$ on $\bar{B}_{p}^{\prime }=\bar{B}(p,r)$ and all the other
values are in the open interval $(0,1)$. We can consider an at most
countable cover of $M_{1}$ with such balls $B_{p}$. In the case when the
cover of $M_{1}$ is a finite set $\{B_{i}\}_{i=\overline{1,r}}$, we can
consider $\varphi =\sum\limits_{i=1}^{r}\psi _{i}$, that is obviously a test
function for $M_{0}$. In the case when the cover of $M_{1}$ is a finite set $%
\{B_{i}\}_{i=\overline{1,r}}$, we can proceed as in \cite[Section 4]{DLPR}.
For each $i\in I\!\!N$ consider the constants $c_{i}$ such that $%
c_{i}\left\Vert \psi _{i}\right\Vert \leq 1/2^{i}$, where the norms are in $%
BC^{\infty }(I\!\!R^{n},I\!\!R)$, then denote $\varphi _{i}=c_{i}\psi _{i}$
and finally 
\begin{equation*}
\varphi =\sum\limits_{i=1}^{\infty }\varphi _{i}.
\end{equation*}%
As in the proof of \cite[Proposition 4.3]{DLPR}, $\varphi $ is a smooth
function and the set of its zeros is $I\!\!R^{n}\backslash M_{1}=M_{0}$, so it is a weak test function for $M_0$.
Using Lemma \ref{lmsmooth} with $\psi _{0}$ given by the formula (\ref{ps0}%
), we obtain a strong test function for $M_{0}$. $\Box $

The existence of a weak test function that is not a strong one depends on
the zero set (i.e. the closed set). For example, the singular foliation of $%
I\!\!R^{n}$ by concentric spheres (as leaves of dimension $n-1$) and the
origin (as a leaf of dimension $0$) is locally regular and the square of the
euclidian norm is a weak test function that is not a strong one.  The
singular foliation having as leaves concentric spheres, as in the previous
example (of dimension $n-1$), outside a compact ball $\bar{B}(\bar{0},\rho
)\subset I\!\!R^{n}$, $\rho >0$, while $\bar{B}(\bar{0},\rho )$ is a union
of points (as leaves of dimension $0$) is also locally regular, but every
weak test function of $\bar{B}(\bar{0},\rho )$ is always a strong one.

\section{The construction of Godbillon-Vey forms and classes}

Integrability conditions for a regular foliation are given by Frobenius
theorem. It can be expressed using differential forms, as, for example, in 
\cite[Ch. 2. and Ch. 3]{To1}. We use this in a similar way as in \cite{Ko}.
If a differentiable $q$-form $\nu $ on $M$ has locally the form $\nu =\omega
_{1}\wedge \cdots \wedge \omega _{q}$, where $\omega _{1}$, $\ldots $, $%
\omega _{q}$ are local one-forms, we say that $\nu $ has rank $q$.

A regular foliation of co-dimension $q$ on a differentiable manifold $M$ is
given by a non-singular global form $\nu \in \Omega ^{q}(M)$ of rank $q$
and, in the locally form $\nu =\omega _{1}\wedge \cdots \wedge \omega _{q}$,
the local one-forms $\omega _{1}$, $\ldots $, $\omega _{q}$ are sections of
the transverse bundle of the foliation, that generate the $\mathcal{F}(M)$%
-module of transverse one-forms (\cite[Proposition 3.9]{To1}). One briefly
say that the foliation (or its tangent bundle) is given by $\nu =0$, or by
vanishing of $\nu $.

Let us consider now two regular foliations $\mathcal{F}_{U}$ and $\mathcal{F}%
_{V}$, $\mathcal{F}_{U|U\cap V}\subset \mathcal{F}_{V|U\cap V}$, such that
the tangent bundles of the foliations $\mathcal{F}_{U}$ and $\mathcal{F}_{V}$
are given of vanishing the differential forms $\omega _{1}\in \Omega
^{q_{1}+q_{2}}(U)$ and $\omega _{2}\in \Omega ^{q_{2}}(V)$ respectively.

\begin{prop}
\label{prGb1}Denoting by $\omega _{1}^{\prime }\in \Omega ^{q_1+q_2}(U\cap V)$
and $\omega_{2}^{\prime }\in \Omega ^{q_2}(U\cap V)$ the
restrictions to $U\cap V$ of $\omega _{1}$ and $\omega _{2}$ respectively,
where $q_{1}>0$, then there is a differentiable form $\theta \in \Omega
^{q_1}(U\cap V)$ such that 
\begin{equation}
\omega _{1}^{\prime }=\omega _{2}^{\prime }\wedge \theta .  \label{Godb1}
\end{equation}
\end{prop}

\emph{Proof.} First, let us suppose that $U=V=U\cap V$ is a domain of
coordinates $\{x^{u},\tilde{x}^{\tilde{u}},\bar{x}^{\bar{u}}\}$, $u=%
\overline{1,p}$, $\tilde{u}=\overline{1,q_{1}}$ and $\bar{u}=\overline{%
1,q_{2}}$ such that $\{x^{u}\}$ and $\{x^{u},\tilde{x}^{\tilde{u}}\}$ are
coordinates on the leaves of $\mathcal{F}_{U|U\cap V}$ and $\mathcal{F}%
_{V|U\cap V}$ respectively. Then $\omega _{1}^{\prime }=h_{1}d\tilde{x}%
^{1}\wedge \cdots \wedge d\tilde{x}^{q_{1}}\wedge d\bar{x}^{1}\wedge \cdots
\wedge d\bar{x}^{q_{2}}$ and $\omega _{2}^{\prime }=h_{2}d\bar{x}^{1}\wedge
\cdots \wedge d\bar{x}^{q_{2}}$ with $h_{1}$, $h_{2}\in \mathcal{F}(U\cap V)$
having no zeros, thus relation (\ref{Godb1}) holds for $\theta =$ $\frac{%
h_{1}}{h_{2}}d\tilde{x}^{1}\wedge \cdots \wedge d\tilde{x}^{q_{1}}$.
Returning to the general case, let us consider a partition of unity $%
\{v_{\alpha }\}_{\alpha \in A}$ on $U\cap V$ subordinated to a cover with
open domain of local foliated charts, as above, where $A$ is finite or $%
A=I\!\!N$. Then define $\theta =\sum\limits_{\alpha \in A}v_{\alpha }\theta
_{\alpha }\in \Omega ^{1}(U\cap V)$. Since $\omega _{1}^{\prime }=\omega
_{2}^{\prime }\wedge \theta _{\alpha }$ and $\sum\limits_{\alpha \in
A}v_{\alpha }=1$, then relation (\ref{Godb1}) holds. $\Box $

In order to avoid coordinates, we consider in the sequel the ideals $%
\mathcal{I}(\mathcal{F}_{U})\subset \Omega ^{\ast }(U)$ and $\mathcal{I}(%
\mathcal{F}_{V})\subset \Omega ^{\ast }(V)$ of differential forms that
vanish when evaluated with all vectors that are tangent to the leaves of $%
\mathcal{F}_{U}$ and $\mathcal{F}_{V}$ respectively. The two ideals are
finitely generated, each homogeneous term containing at least one of the
local forms that on $U\cap V$ can be taken of the form $\{\tilde{\omega}^{\tilde{u}},\bar{\omega}^{\bar{%
u}}\}_{\tilde{u}=\overline{1,q_{1}},\bar{u}=%
\overline{1,q_{2}}}$ and $\{\bar{\omega}^{\bar{u}}\}_{\bar{u}=\overline{%
1,q_2}}$ respectively. Notice that $d\bar{\omega}^{\bar{u}}=\sum\limits_{%
\bar{v}=1}^{q_{2}}\bar{\omega}^{\bar{v}}\wedge \nu _{\bar{v}}^{\bar{u}}$ and 
$d\tilde{\omega}^{\tilde{u}}=\sum\limits_{\bar{v}=1}^{q_{2}}\bar{\omega}^{%
\bar{v}}\wedge \nu _{\bar{v}}^{\tilde{u}}+\sum\limits_{\tilde{v}=1}^{q_{1}}%
\tilde{\omega}^{\tilde{v}}\wedge \nu _{\tilde{v}}^{\tilde{u}}$, with $\nu _{%
\bar{v}}^{\bar{u}}$, $\nu _{\bar{v}}^{\tilde{u}}$ and $\nu _{\tilde{v}}^{%
\tilde{u}}\in \Omega ^{1}(U\cap V)$. Then $\omega _{2}$ has the local form 
\begin{equation}
\omega _{2}=h_{2}\bar{\omega}^{1}\wedge \cdots \wedge \bar{\omega}^{q_2}.
\label{formom2}
\end{equation}

The Frobenius theorem used for $\mathcal{F}_{U}$ and $\mathcal{F}_{V}$ reads
that there are $\mu _{1}\in \Omega ^{1}(U)$ and $\mu _{2}\in \Omega ^{1}(V)$
such that 
\begin{equation}
d\omega _{1}=\omega _{1}\wedge \mu _{1},\ d\omega _{2}=\omega _{2}\wedge \mu
_{2}.  \label{Godb2}
\end{equation}

A product of $q_{1}+q_{2}+1$ forms in $\mathcal{I}(\mathcal{F}_{U})$ as well
as of $q_{2}+1$ forms in $\mathcal{I}(\mathcal{F}_{V})$ are null. This
enables to consider the closed \emph{Godbillon-Vey forms }$\mu _{1}\wedge
(d\mu _{1})^{q_{1}+q_{2}}\in \Omega ^{2\left( q_{1}+q_{2}\right) +1}(U)$ and 
$\mu _{2}\wedge (d\mu _{2})^{q_{2}}\in \Omega ^{2q_{2}+1}(V)$ and the \emph{%
Godbillon-Vey classes} of the foliations $\mathcal{F}_{U}$ and $\mathcal{F}%
_{V}$ as the cohomology classes $[\mu _{1}\wedge (d\mu
_{1})^{q_{1}+q_{2}}]\in H^{2\left( q_{1}+q_{2}\right) +1}(U)$ and $[\mu
_{2}\wedge (d\mu _{2})^{q_{2}}]\in H^{2q_{2}+1}(V)$.

Let us look closely to $U\cap V$, when the relation (\ref{Godb1}) holds. For
sake of simplicity, we use notations $\omega _{1}$ and $\omega _{2}$ instead
of $\omega _{1}^{\prime }$ and $\omega _{2}^{\prime }$ respectively.

Differentiating by $d$ (\ref{Godb1}), then using (\ref{Godb2}) and the usual
properties of the exterior product, we obtain%
\begin{equation*}
\omega _{2}\wedge \left( \left( -1\right) ^{q_{2}}d\theta -\theta \wedge
\left( \mu _{1}-\left( -1\right) ^{q_1}\mu _{2}\right) \right) =0.
\end{equation*}%
Taking into account (\ref{formom2}), then 
\begin{equation}
d\theta -\left( -1\right) ^{q_{2}}\theta \wedge \left( \mu _{1}-\left(
-1\right) ^{q_1}\mu _{2}\right) =\sum\limits_{\bar{v}=1}^{q_{2}}\bar{%
\omega}^{\bar{v}}\wedge \eta _{\bar{v}},  \label{Godb3}
\end{equation}
with $\eta _{\bar{v}}\in \Omega ^{q_{1}}(U\cap V)$. Thus the left side of
equality (\ref{Godb3}) belongs to $\mathcal{I}(\mathcal{F}_{V})_{|U\cap
V}\subset \Omega ^{\ast }(U\cap V)$. Denote by 
\begin{equation}
\mu _{3}=\left( -1\right) ^{q_{2}}\left( \mu _{1}-\left( -1\right)
^{q_1}\mu _{2}\right) .  \label{mu3}
\end{equation}

Differentiating by $d$ and using again the same relation (\ref{Godb3}), we
obtain 
\begin{equation}
\theta \wedge d\mu _{3}=\sum\limits_{\bar{v}=1}^{q_{2}}\bar{\omega}^{\bar{v}%
}\wedge \bar{\eta}_{\bar{v}},  \label{Godb4}
\end{equation}%
with $\bar{\eta}_{\bar{v}}\in \Omega ^{q_{1}+1}(U\cap V)$, i.e. $\theta
\wedge d\mu _{3}\in \mathcal{I}(\mathcal{F}_{V})_{|U\cap V}$. But using
local coordinates as in the proof of Proposition \ref{prGb1}, we have that,
on a domain $U^{\prime }$ of such coordinates, there is a local function $%
h_{3}$ such that $\theta -h_{3}d\tilde{x}^{1}\wedge \cdots \wedge d\tilde{x}%
^{q_{1}}\in \mathcal{I}(\mathcal{F}_{V})_{|U^{\prime }}$. Using this fact in
(\ref{Godb4}),  for $U^{\prime}=U\cap V$, it follows that 
\begin{equation*}
d\mu _{3}=\sum\limits_{\bar{v}=1}^{q_{2}}\bar{\omega}^{\bar{v}}\wedge \tilde{%
\eta}_{\bar{v}}
\end{equation*}%
with $\tilde{\eta}_{\bar{v}}\in \Omega ^{1}(U\cap V)$, i.e. $d\mu _{3}\in 
\mathcal{I}(\mathcal{F}_{V})_{|U\cap V}$. But $d\mu _{2}\in \mathcal{I}(%
\mathcal{F}_{V})_{|U\cap V}$, thus using (\ref{mu3}) it follows that $d\mu
_{1}\in \mathcal{I}(\mathcal{F}_{V})_{|U\cap V}$.

\begin{prop}
\label{prgb01}Assuming $q_{1}>0$, then the following assertions hold true:
 
\begin{enumerate}

\item The Godbillon-Vey form  and its cohomology class
according to the foliation $\mathcal{F}_{U|U\cap V}$, both vanish.

\item If $\mathcal{F}^{\prime }\subset \mathcal{F}^{\prime \prime }$, $%
\mathcal{F}^{\prime }\neq \mathcal{F}^{\prime \prime }$, are regular
foliations on $M$ and the foliation $\mathcal{F}^{\prime \prime }$ has not a
null co-dimension, then the Godbillon-Vey class of $\mathcal{F}^{\prime }$
vanishes.
\end{enumerate}
\end{prop}

\emph{Proof.}  If $q_{1}>0$, then $%
q_{1}+q_{2}\geq q_2+1$, thus $\left( d\mu _{1}\right) ^{q_{1}+q_{2}}=0$ because $\left( d\mu _{1}\right) ^{1+q_{2}}=0$; it follows that $\mu
_{1}\wedge \left( d\mu _{1}\right) ^{q_{1}+q_{2}}=0$, as well as its
cohomology class, thus 1. follows. Then 2. is a simple consequence of 1. $
\Box $

The result in  Proposition \ref{prgb01}  allows to consider the Godbillon-Vey class of
the foliation $\mathcal{F}_{U_{1}}$ having the maximal co-dimension $q_{\max
}=m-r_{\min }$, on the open subset $U_{r_{\min }}\subset M$; the foliation
has the leaves of minimal dimension. The Godbillon-Vey class is the class $%
[\mu _{r_{\min }}\wedge \left( d\mu _{r_{\min }}\right) ^{m-r_{\min }}]$. The
differential form $GV_{r_{\min }}=\mu _{r_{\min }}\wedge \left( d\mu
_{r_{\min }}\right) ^{m-r_{\min }}\in \Omega ^{1+2q_{\max }}(U_{r_{\min} })$ is
null on any intersection $U_{r_{\min }}\cap U_{0}\neq \emptyset $, where $%
U_{0}$ is an open subset corresponding to a foliation $\mathcal{F}_{U_{0}}$
of co-dimension $q_{0}=m-r_{0}<q_{\max }=m-r_{\min }$. Thus, extending $%
GV_{r_{\min }}$ as null outside $U_{r_{\min }}$, we obtain a global closed
form that gives $GV_{\min }(\mathcal{F})\in H^{1+2q_{\max }}(M)$; we call it
as the \emph{Godbillon-Vey class on leaves of minimal dimension} of the
locally regular foliation $\mathcal{F}$. 

In the general case, let us consider the ascending sequence of open sets $%
\Sigma _{\geq r_{k}}\subset \Sigma _{\geq r_{k-1}}\subset \cdots \subset
\Sigma _{\geq r_{1}}\subset \Sigma _{\geq r_{0}}=M$. Denote by $\mathcal{F}%
_{\Sigma _{\geq r_{i}}}$ the restriction of $\mathcal{F}$ to the open set $%
\Sigma _{\geq r_{i}}$, $i=\overline{0,k}$; notice that the set $\Sigma
_{\geq r_{i}}$ is saturated by the leaves of $\mathcal{F}=\mathcal{F}%
_{\Sigma _{\geq r_{0}}}$. The subset $\Sigma _{r_{i}}\subset \Sigma _{\geq
r_{i}}$ is that of minimal dimensions of leaves. We can consider the
Godbillon-Vey classes $GV_{\min }(\mathcal{F}_{\Sigma _{\geq r_{i}}})\in
H^{2\left( m-r_{i}\right) +1}(\Sigma _{\geq r_{i}})$. In particular, 
\begin{displaymath}
GV_{\min }(\mathcal{F})=GV_{\min }(\mathcal{F}_{\Sigma _{\geq r_{0}}})\in
H^{2\left( m-r_{0}\right) +1}(\Sigma _{\geq r_{0}})=H^{2\left(
m-r_{0}\right) +1}(M).
\end{displaymath}

\begin{theor}
\label{thmain}A Godbillon-Vey form of the leaves extends to a global
cohomology class $GV_{\min }(\mathcal{F})\in H^{1+2q_{\max }}(M)$ (for the
leaves of minimal dimension) and to some Godbillon-Vey classes $GV_{\min }(%
\mathcal{F}_{\Sigma _{\geq r_{i}}})\in H^{2\left( m-r_{i}\right) +1}(\Sigma
_{\geq r_{i}})$ (for the leaves on the other $U_{i}$, $i>0$).
\end{theor}

In order to obtain global cohomology classes on $M$, the construction on the
Godbillon-Vey class on the leaves of minimal dimension can be extended to
the other strata, provided that there is a foliated test function according
to that stratum. We perform below this construction.

Let us suppose that the foliation $\mathcal{F}_{r_{i}}$ on $U_{i}\subset M$
has the dimension $r_{i}$ of leaves and it is defined on $U_{i}$ by the
equation $\omega _{i}=0$, where $\omega _{i}\in \Omega ^{q_{i}}(U_{i})$, $%
q_{i}=m-r_{i}$. Then%
\begin{equation*}
d\omega _{i}=\omega _{i}\wedge \mu _{i}
\end{equation*}%
with $\mu _{i}\in \Omega ^{1}(U_{i})$. We suppose below that there is a test
function $\varphi _{i}\in \mathcal{F}(M)$ for $M\backslash U_{i}$ that
restricts to a basic function for the foliation $\mathcal{F}_{r_{i}}$ on $%
U_{i}$; we suppose also that $\bar{\mu}_{i}=\varphi _{i}\mu _{i}$ (where $%
\mu _{i}$ is defined by zero on $M\backslash U_{i}$) is differentiable on $M$%
, i.e. $\bar{\mu}_{i}\in \Omega ^{1}(M)$; this is always true if $\varphi
_{i}$ is a strong test function.

\begin{prop}
\label{prbaz}Let us suppose that the test function $\varphi _{i}$ is basic
and $\bar{\mu}_{i}=\varphi _{i}\mu _{i}$ is differentiable on $M$. Then the
differential form $\bar{\nu}_{i}=\bar{\mu}_{i}\wedge \left( d\bar{\mu}%
_{i}\right) ^{q_{i}}$ is closed, giving a cohomology class $[\bar{\nu}%
_{i}]\in H^{2q_{i}+1}(M)$.
\end{prop}

\emph{Proof.} We have $\bar{\nu}_{i}=\bar{\mu}_{i}\wedge \left( d\bar{\mu}%
_{i}\right) ^{q_{i}}=\varphi _{i}^{1+q_{i}}\mu _{i}\wedge \left( d\mu
_{i}\right) ^{q_{i}}$.  But, if $\varphi _{i}$ is basic, then $\psi _{i}=\varphi
_{i}^{1+q_{i}}$ is also basic and $d\psi _{i}\wedge \mu _{i}\wedge \left(
d\mu _{i}\right) ^{q_{i}}=0$. Thus $d\bar{\nu}_{i}=0$ and the conclusion
follows. $\Box $

Notice that if the maximal stratum has the dimension $r_{k}=m$, then its
Godbillon-Vey form vanishes, as well as its Godbillon-Vey class. In
particular, if a family of regular foliations has $\mathcal{R}%
=\{r_{0},r_{1}\}$ and $r_{1}=m$, then the only possible non-null is the
Godbillon-Vey class of the leaves of minimal dimension.

\section{Two cases}

First, we prove that the usual Godbillon-Vey class of a regular foliation is
the same with the Godbillon-Vey class of leaves of minimal dimension of a
suitable non-trivial family\ of regular foliations. Let $(M,\mathcal{F}_{0})$
be a regular foliation of co-dimension $q_{0}$ defined by a $q_{0}$%
--differential form $\omega _{0}=0$, such that $d\omega _{0}=\omega
_{0}\wedge \mu _{0}$. Let us consider two open and non-void subsets $W$, $%
U_{2}$ having the properties that $\bar{W}\subset U_{2}$ and $\varphi \in 
\mathcal{F}(M)$ is a Uryson function such that $supp$ $\varphi =M\backslash 
\bar{W}=U_{1}$. Consider on $U_{1}$ the foliation $\mathcal{F}_{U_{1}}$ as
being the restriction to $U_{1}$ of foliation $\mathcal{F}_0$. Let us suppose
that there is on $U_{2}$ a non--trivial foliation $\mathcal{F}_{U_{2}}$ such
that its leaves are saturated by leaves of $\mathcal{F}_{0|U_{2}}$ (for this
we can take $U_{2}$ the domain of a $\mathcal{F}_{0}$--foliate simple chart
and then take as $\mathcal{F}_{U_{2}}$ a proper foliation having as
subfoliation $\mathcal{F}_{0|U_{2}}$, for example, a trivial foliation with
one leaf). The foliation $\mathcal{F}_{U_{2}}$ is defined by the $q_{0}$%
--form $\tilde{\omega}=\varphi \omega _{0}$, that has the same support as $%
\varphi $. The foliations $\mathcal{F}_{U_{1}}$ and $\mathcal{F}_{U_{2}}$
give a non-trivial family\ of regular foliations on $M$. The Godbillon-Vey
class $GV_{\min }(\mathcal{F})\in H^{2q_{0}+1}(M)$ is given extending
naturally (using Proposition \ref{prgb01}) a form that gives the
Godbillon-Vey class of $\mathcal{F}_{U_{1}}$.

\begin{prop}
\label{prex01}The Godbillon-Vey class $GV_{\min }(\mathcal{F})$ is the same
as $GV(\mathcal{F}_{0})$, the usual Godbillon-Vey class of $\mathcal{F}_{0}$.
\end{prop}

\emph{Proof.} The Godbillon-Vey class of $\mathcal{F}_{0}$ is given by $\left[\mu_0\wedge(d\mu_0)^{q_0}\right]$, where the definition does not depend of $\omega_0 $ and $\mu_0$ (see \cite[Theorem 3.11]{To1}). It can be easy proved that we
can take the restriction of $\omega_0 $ to $U_{2}$ having the form  $fd%
\bar{x}^{1}\wedge \cdots \wedge d\bar{x}^{q_0}$, where $\{\bar{x}^{\bar{%
u}}\}_{\bar{u}=\overline{1,q_{0}}}$ are transverse coordinates for $\mathcal{
F}_{0}$ on $U_{2}$, thus $\mu_0|_{U_{2}}=(-1)^{q_{0}}d\log f$  and $d\mu_0|
_{U_{2}}=0$. Thus the restriction of the differential form $\mu_0 \wedge
(d\mu_0 )^{q_{0}}$ to $U_{2}$ vanishes and it extends the differential form
on $U_{1}$ that gives the Godbillon-Vey class of $\mathcal{F}_{0|U_{2}}$,
thus it gives $GV_{\min }(\mathcal{F})$. It follows that $GV_{\min }(%
\mathcal{F})=GV(\mathcal{F}_{0})$. $\Box $

We consider below a non-trivial case when the Godbillon-Vey class vanishes.
More specifically, we prove that for a regular (weak) test function $%
\varphi _{i}\in \mathcal{F}(M)$ for $M\setminus U_{i}$ that restricts to a
basic function for the foliation $\mathcal{F}_{r_{i}}$ on $U_{i}$ the
cohomology class $[\overline{\nu }_{i}]\in H^{2q_{i}+1}(M)$ vanishes. 

Firstly we shall need some preliminary notions about singular forms and
cohomology attached to a function, for more see \cite{Mo1, Mo2}.
Accordingly, for a smooth function $f\in \mathcal{F}(M)$ and $U\subset M$ a $%
p$--form $\omega \in \Omega ^{p}(U)$ is called a \textit{singular }$p$%
\textit{--form }if the form $f^{p}\omega $ can be extended to a smooth form
on $M$, that is $f^{p}\omega \in \Omega ^{p}(M)$. We denote the space of
singular $p$--forms with respect to $f$ by $\Omega _{f}^{p}(M)$. We notice
that if $\omega \in \Omega _{f}^{p}(M)$ then $d\omega \in \Omega
_{f}^{p+1}(M)$ and so we have a differential complex $\left( \Omega
_{f}^{\bullet }(M),d\right) $. The cohomology of this differential complex
is isomorphic with the cohomology attached to the function $f$, denoted by $%
H_{f}^{\bullet }(M)$, which is defined as cohomology of the differential
complex $\left( \Omega ^{\bullet }(M),d_{f}\right) $, where the coboundary
operator $d_{f}:\Omega ^{p}(M)\rightarrow \Omega ^{p+1}(M)$ is defined by $%
d_{f}\omega =fd\omega -pdf\wedge \omega $. The mentioned isomorphism is
produced by the map of chain complexes $\phi :\left( \Omega _{f}^{\bullet
}(M),d\right) \rightarrow \left( \Omega ^{\bullet }(M),d_{f}\right) $ given
by $\phi ^{p}:\Omega _{f}^{p}(M)\rightarrow \Omega ^{p}(M)$, $\phi (\omega
)=f^{p}\omega $, see \cite{Mo2}.

Now, let us return to our study. As well as we seen from the above
discussion $\mu_i\in\Omega^1_{\varphi_i}(M)$ and, accordingly $%
d\mu_i\in\Omega^2_{\varphi_i}(M)$. We have then that $\mu_i\wedge
(d\mu_i)^{q_i}\in\Omega^{2q_i+1}_{\varphi_i}(M)$. Since $\mu_i\wedge
(d\mu_i)^{q_i}$ is closed, from the above isomorphism we have that $%
\varphi_i^{2q_i+1}\mu_i\wedge(d\mu_i)^{q_i}$ is $d_{\varphi_i}$--closed.
Thus, if $\varphi_i$ is basic function for the foliation $\mathcal{F}%
_{r_{i}} $ on $U_{i}$ then $d_{\varphi_i}(\varphi_i^{q_i}\overline{\nu}_i)=0$
which leads to the cohomology class $[\varphi_i^{q_i}\overline{\nu}_i]\in
H^{2q_i+1}_{\varphi_i}(M)$. Let us consider now the \textit{regular} case
for the test function $\varphi_i$, that is $\varphi_i$ does not have
singularities in a neighborhood of its zero set (i.e., $0$ is a regular
value). The subsets $S_i=\varphi_i^{-1}(\{0\})=M\setminus U_i$ are then
embedded submanifolds of $M$. We also assume that $S_i$ are connected.

We consider some useful notations. Let $V_i\subset V^{\prime}_i$ be tubular
neighborhoods of $S_i$. We may assume that $V_i=S_i\times
]-\varepsilon_i,\varepsilon_i[$ and $V^{\prime}_i=S_i\times]-\varepsilon^{%
\prime}_i,\varepsilon^{\prime}_i[$, with $\varepsilon^{\prime}_i>%
\varepsilon_i$, and that $\varphi_i|_{V^{\prime}_i}:S_i\times]-\varepsilon^{%
\prime}_i,\varepsilon^{\prime}_i[\rightarrow \mathbb{R}$, $(x,t)\mapsto t$.
We denote by $\pi_i$ the projections $V^{\prime}_i\rightarrow S_i$. Let $%
\rho:\mathbb{R}\rightarrow\mathbb{R}$ be a smooth function which is $1$ on $%
[-\varepsilon_i,\varepsilon_i]$ and has support contained in $%
[-\varepsilon^{\prime}_i,\varepsilon^{\prime}_i]$. Note that the function $%
\rho\circ \varphi_i$ is $1$ on $V_i$, and we can assume that the function $%
\rho\circ\varphi_i$ vanishes on $M\setminus V^{\prime}_i$. If $\omega$ is a
form on $S_i$, we will denote by $\widetilde{\omega}$ the form $%
\rho(\varphi_i)\pi_i^*\omega$ and notice that $d\varphi_i\wedge d\widetilde{%
\omega}=d\varphi_i\wedge\widetilde{d\omega}$, see \cite{Mo2}.

According to Theorem 4.1 from \cite{Mo2}, if $0$ is a regular value of $%
\varphi_i$ then, for each $p\geq 1$, there is an isomorphism 
\begin{equation}  \label{iso}
H^p_{\varphi_i}(M)\cong H^p(M)\oplus H^{p-1}(S_i),
\end{equation}
given by $\Phi:\Omega^p(M)\oplus\Omega^{p-1}(S_i)\rightarrow\Omega^p(M)$
defined by $\Phi(\alpha,\beta)=\varphi_i^p\alpha+\varphi_i^{p-1}d\varphi_i%
\wedge\widetilde{\beta}$.

Now, taking into account the isomorphism \eqref{iso} it follows that there
exist $\alpha _{i}\in \Omega ^{2q_{i}+1}(M)$ and $\beta _{i}\in \Omega
^{2q_{i}}(S_{i})$ with $d\alpha _{i}=d\beta _{i}=0$ such that 
\begin{equation}
\varphi _{i}^{q_{i}}\overline{\nu }_{i}=\varphi _{i}^{1+2q_{i}}\alpha
_{i}+\varphi _{i}^{2q_{i}}d\varphi _{i}\wedge \widetilde{\beta }_{i}.
\label{I1}
\end{equation}%
Thus we obtain that $\alpha _{i}=\varphi _{i}^{-1-q_{i}}\overline{\nu }_{i}-%
\frac{d\varphi _{i}}{\varphi _{i}}\wedge \widetilde{\beta }_{i}$ and by
differentiation and taking into account $d\overline{\nu }_{i}=d\alpha
_{i}=d\beta _{i}=0$, one get 
\begin{equation*}
(-1-q_{i})\varphi _{i}^{-2-q_{i}}d\varphi _{i}\wedge \overline{\nu }_{i}=0,
\end{equation*}%
where we have used $d\varphi _{i}\wedge d\widetilde{\beta }_{i}=d\varphi
_{i}\wedge \widetilde{d\beta _{i}}=0$.

Now, since $d\overline{\nu }_{i}=0$ and $d\varphi _{i}\wedge \overline{\nu }%
_{i}=0$, by Proposition 3.4 from \cite{Mo1} there exist $\overline{\tau }%
_{i}\in \Omega ^{2q_{i}-1}(M)$ such that $\overline{\nu }_{i}=d\varphi
_{i}\wedge d\overline{\tau }_{i}$ and so $\overline{\nu }_{i}=d(\varphi _{i}d%
\overline{\tau }_{i})$. Thus, we obtain the announced result:

\begin{prop}
\label{prex02}If $0$ is a regular value for the (weak) test function $%
\varphi _{i}$ that is also basic, then the cohomology class $[\overline{\nu }%
_{i}]\in H^{2q_{i}+1}(M)$ vanishes.
\end{prop}

\section*{Acknowledgement}
The author CI is supported by the Sectorial Operational Program Human Resources Development (SOP HRD), financed from the European Social Fund and by the Romanian Government under the Project number POSDRU/159/1.5/S/134378.

\noindent
Cristian Ida\\
Department of Mathematics and Computer Science,\\
University Transilvania of Bra\c{s}ov, \\
Faculty of Mathematics and Computer Science,\\
Address: Bra\c{s}ov 500091, Str. Iuliu Maniu 50, Rom\^{a}nia.\\
E-mail:\textit{cristian.ida@unitbv.ro}
\[
\]
\noindent 
Paul Popescu\\
Department of Applied Mathematics,\\
University of Craiova, Faculty of Exact Sciences,\\
13, A.I.Cuza st., Craiova, 200585, Rom\^{a}nia.\\
E-mail:\textit{paul$_{-}$p$_{-}$popescu@yahoo.com}


\begin{thebibliography}{99}
\bibitem{Da} Dazord P., \emph{F\'{e}uilletages \`{a} singularit\'{e}s},
Indagationes Math. Volumen, 1, 47 (1985) 21-39.

\bibitem{Do} Dominguez, D., \emph{Sur les Classes Caract\`{e}ristiques des
Sous-Feuilletages}, Publ. RIMS, Kyoto Univ., 23 (1987) 813-840.


\bibitem{DLPR} Drager L.D., Lee J.M., Park E., Richardson K., \emph{Smooth
vector subbundles are finitely generated}, Ann. Glob. Anal. Geom., 41, 3
(2012) 357-369.

\bibitem{GHV} Greub W., Halperin S., Vanstone R., \emph{Connections,
Curvature, and Cohomology}, vol.I, Academic Press, New York, 1972.

\bibitem{Go} Godbillon C, Reeb G., \emph{Feuilletages: \'{e}tudes g\'{e}om%
\'{e}triques}, Birkh\"{a}user, Basel, 1991.

\bibitem{Hi} Hirsch M., \emph{Differential Topology}, Graduate Text in Math.
33 Springer-Verlag, New York, 1976.

\bibitem{Ho} Hoster, M. \emph{Derived secondary classes for flags of
foliations}, PhD Diss. LMU, 2001.

\bibitem{Lee} Lee J.M., \emph{Introduction to Smooth Manifolds}, Springer
Verlag, New York, 2003.

\bibitem{Mo1} Monnier Ph., \emph{Computations of Nambu-Poisson cohomologies}%
. IJMMS 26, 2 (2001) 65--81.

\bibitem{Mo2} Monnier Ph., \emph{A cohomology attached to a function}.
Diff. Geom. and Appl. 22 (2005) 49--68.

\bibitem{Mu} Muger M., \emph{An Introduction to Differential Topology, de
Rham Theory and Morse Theory}, http://www.math.ru.nl/\symbol{126}%
mueger/diff\_notes.pdf (2005).

\bibitem{Ko} Kotschick D., \emph{Godbillon-Vey invariants for families of
foliations}, Eliashberg, Yakov (ed.) et al., Symplectic and contact
topology: Interactions and perspectives. Papers of the workshop on
symplectic and contact topology, quantum cohomology, and symplectic field
theory, Montreal and Toronto, Canada, March--April 2001. Providence, RI,
AMS, Fields Inst. Commun. 35 (2003) 131-144.

\bibitem{La} Lang S., \emph{Differential and Riemannian Manifolds,} 3-rd
ed., Springer Verlag, New York, 1995.

\bibitem{To1} Tondeur P., \emph{Foliations on Riemannian manifolds},
Universitext. Springer-Verlag, Berlin, Heidelberg, New York, 1988.

\bibitem{Va1} Vaisman I., \emph{Lectures on the Geometry of Poisson Manifolds%
}, Progress in Math., vol. 118, Birkh\"{a}user Verlag, Boston, 1994.

\bibitem{Wi} H. Whitney, \emph{Analytic extensions of differentiable
functions defined in closed sets}, Trans. Amer. Math. Soc. 36 (1934) 63-89.
\end{thebibliography}
\end{document}